\documentclass[a4paper]{scrartcl}
\usepackage{amsmath,amssymb,url,graphicx, verbatim}
\usepackage{varioref}
\usepackage{relsize}
\newtheorem{theorem}{Theorem}[section]
\newtheorem{lemma}[theorem]{Lemma}

\DeclareMathOperator{\Tr}{Tr}
\numberwithin{equation}{section}

\begin{document}

\begin{center}
 \LARGE \textbf{SPECTRAL DISTRIBUTION OF NON-INDEPENDENT RANDOM MATRIX ENSEMBLES INDUCED BY LACUNARY SYSTEMS}
\end{center}

\vspace{3ex}

\begin{center}
 THOMAS L\"OBBE
\footnote{The results are part of the author's PhD thesis supported by IRTG 1132, University of Bielefeld}
\end{center}

\vspace{3ex}

\begin{abstract}
ABSTRACT. For two lacunary sequences $(M_{n,1})_{n\geq 2},(M_{n,2})_{n\geq 0}$ and suitable functions $f$ we introduce random matrix ensembles with 
\begin{equation*}
 X_{n,n'}=f(M_{n+n',1}x_1,M_{|n-n'|,2}x_2).
\end{equation*}
We prove weak convergence of the mean empirical eigenvalue distribution towards the semicircle law under some further number theoretic properties of the sequence $(M_{n,1})_{n\geq 1}$. To prove this result we show (\ref{252}). Furthermore we give examples to show that even in this particular class of random matrix ensembles the asymptotic behaviour of the spectrum becomes delicate. We prove that the empirical spectral distribution does not converge to the semicircle law in general even if the correlation of two entries decays exponentially in the distance. For $f(x_1,x_2)=1/\sqrt{2}\cdot(\cos(2\pi(x_1+x_2))+\cos(4\pi(x_1+x_2)))$ and $M_{n,1}=2^n$ we show that the mean empirical spectral distribution does not converge to semicircle law while for any sequence $(M_{n,1})_{n\geq 1}$ with $M_{n+1,1}/M_{n,1}\to\infty$ for $n\to\infty$ and any periodic function $f$ of finite total variation in the sense 
of Hardy and Krause with mean zero and unit variance the mean spectral distribution converges to the semicircle law. 
\end{abstract}

\vspace{5ex}

\section{Introduction}

\subsection*{Lacunary sequences}
Let $(M_n)_{n\geq 1}$ be a sequence of integers satisfying a Hadamard gap type condition of the form
\begin{equation}
\label{17}
M_{n+1}\geq qM_n
\end{equation}
for all $n\in\mathbb{N}$ and some absolute constant $q>1$.

For such a sequence Salem and Zygmund \cite{SZ47} proved that for any sequence of integers $(a_n)_{n\geq 1}$ with
\begin{equation*}
 a_N=o(A_N) \quad \textnormal{for} \quad A_N=\frac{1}{2}\left(\sum_{n=1}^Na_n^2\right)^{1/2}
\end{equation*}
we have
\begin{equation}
\lim_{N\to\infty}\mathbb{P}\left(\frac{1}{A_N}\sum_{n=1}^Na_n\cos(2\pi M_nx)\leq t \right)=\Phi(t)
\end{equation}
where $\mathbb{P}$ denotes the probability measure induced by the Lebesgue measure on $[0,1)^d$ and $\Phi$ denotes the standard normal distribution, i.e. for all $t\in\mathbb{R}$ we have
\begin{equation*}
 \Phi(t)=\frac{1}{\sqrt{2\pi}}\int_{-\infty}^{t}e^{-\frac{1}{2}y^2}\,dy.
\end{equation*}

Furthermore Weiss \cite{W59} (see also Salem and Zygmund \cite{SZ50}, Erd\H{o}s and G\'al \cite{EG55}) showed that
\begin{equation}
\limsup_{N\to\infty}\frac{\sum_{n=1}^Na_n\cos(2\pi M_nx)}{\sqrt{2A_N^2\log\log(N)}}=1 \quad a.e.
\end{equation}
under the condition
\begin{equation*}
 a_N=o\left(\frac{A_N}{\sqrt{\log(\log(A_N))}}\right).
\end{equation*}
With coefficients satisfying $a_N=o(A_N^{1-\delta})$ for some $\delta>0$ Philipp and Stout \cite{PS75} showed that there exists a Brownian Motion $\{W(t):t\geq 0\}$ such that
\begin{equation*}
 \sum_{n=1}^Na_n\cos(2\pi M_nx)=W(A_N)+\mathcal{O}(A_N^{1/2-\varrho}) \quad \textnormal{a.e.}
\end{equation*}
for some $\varrho>0$.
Therefore for lacunary $(M_n)_{n\geq 1}$ the sequence $(a_n\cos(2\pi M_nx))_{n\geq 1}$ shows a behaviour typical for independent, identically distributed random variables. One could ask whether this holds for other periodic
functions as well. The answer is negative in general. By a result of Erd\H{o}s and Fortet (see \cite{K49}) for $f(x)=\cos(2\pi x)+\cos(4\pi x)$ and $M_n=2^n-1$ we have
\begin{equation}
\lim_{N\to\infty}\mathbb{P}\left(\frac{1}{\sqrt{N}}\sum_{n=1}^Nf(M_nx)\leq t\right)=\frac{1}{\sqrt{\pi}}\int_{0}^1\int_{-\infty}^{t|\cos(\pi s)|/2}e^{-u^2}\,duds
\end{equation}
and
\begin{equation}
\limsup_{N\to\infty}\frac{\sum_{n=1}^Nf(M_nx)}{\sqrt{N\log(\log(N))}}=2\cos(\pi x) \quad a.e.
\end{equation}
Thus neither the Central Limit Theorem nor the Law of the Iterated Logarithm is satisfied. This result was later generalized by Conze and Le Borgne \cite{CB11} (see also \cite{A13a} for further information).
On the other hand Kac \cite{K46} showed that any one-periodic function $f:\mathbb{R}\to\mathbb{R}$ of mean zero which is of bounded variation on $[0,1)$ or Lipschitz-continuous satisfies
\begin{equation}
\lim_{N\to\infty}\mathbb{P}\left(\frac{1}{\sqrt{N}}\sum_{n=1}^Nf(2^nx)\leq t\sigma\right)=\Phi(t)
\end{equation}
if
\begin{equation}
\sigma^2=\mathbb{E}[f]+2\sum_{n=1}^{\infty}\mathbb{E}[f(x)f(2^nx)]\neq 0.
\end{equation}
Furthermore Maruyama \cite{M50} and Izumi \cite{I51} proved
\begin{equation}
 \limsup_{N\to\infty}\frac{\sum_{n=1}^Nf(2^nx)}{\sqrt{2N\log(\log(N))}}=\sigma \quad \textnormal{a.e.}
\end{equation}
This illustrates that the behaviour of $(f(M_nx))_{n\geq 1}$ does not only depend on the speed of growth of $(M_n)_{n\geq 1}$ but also on number theoretic properties of the sequence $(M_n)_{n\geq 1}$.
Later on the Central Limit Theorem was shown for more general lacunary sequences. By a result of Gaposhkin \cite{G66}
\begin{equation}
 \label{1}
 \lim_{N\to\infty}\mathbb{P}\left(\sum_{n=1}^Nf(M_nx)\leq t\sigma_N\right)=\Phi(t)
\end{equation}
holds for sequences $(M_n)_{n\geq 1}$ satisfying
\begin{equation}
 \label{2}
 \sigma_N^2=\mathlarger{\mathlarger{\int}}_0^1\left(\sum_{i=1}^Nf(M_nx)\right)^2\,dx\geq CN,
\end{equation}
for an absolute constant $C>0$ and one of the following conditions
\begin{itemize}
 \item $\frac{M_{n+1}}{M_n}\in\mathbb{N}$, \quad \textnormal{for all } $n\in\mathbb{N}$,
 \item $\lim_{n\to\infty}\frac{M_{n+1}}{M_n}=\theta$, \quad such that $\theta^r$ irrational for all $r\in\mathbb{N}$.
\end{itemize}
Takahashi \cite{T61} showed (\ref{1}) for $M_{n+1}/M_{n}\to\infty$ and $\alpha$-Lipschitz-continuous functions.
The connection between the Central Limit Theorem and the number of solutions of certain Diophantine equations is due to Gaposhkin \cite{G70}. 
Consider the linear Diophantine equation
\begin{equation*}
 aj\pm a'j'=\nu
\end{equation*}
for fixed integers $j,j',\nu$. In general the set of solutions consists of all pairs of integers $a,a'$ such that equality holds but we restrict ourselves to those solutions with $a=M_n$ and $a'=M_{n'}$ for $n,n'\in\mathbb{N}$ and rather regard the indices $n,n'$ as solutions of this equation.
The Central Limit Theorem holds for lacunary sequences $(M_n)_{n\geq 1}$ satisfying (\ref{2}), if for any fixed $j,j',\nu$ the number of solutions of the Diophantine equation
\begin{equation}
\label{61}
 M_nj\pm M_{n'}j'=\nu
\end{equation}
is bounded by an absolute constant $C_{j,j'}>0$ which is independent of $\nu$. Observe that ``nice'' periodic functions can be approximated by trigonometric polynomials very well. Thus because of the product-to-sum identities of trigonometric functions the behaviour of the moments of $\sum f(M_nx)$ depends on the number of solutions of Diophantine equations of certain length.
Aistleitner and Berkes \cite{AB10} improved this result: For a lacunary sequence $(M_n)_{n\geq 1}$ satisfying the Hadamard gap condition set
\begin{equation}
 \label{22}
 \begin{aligned}
  L(N,G,\nu) \:\:\:  = & \:\:\: |\{1\leq n,n'\leq N:\\
  &\:\:\: \exists j,j'\in\mathbb{Z}^d,1\leq ||j||_{\infty}|,||j'||_{\infty}\leq G,M_n^Tj\pm M_{n'}^Tj'=\nu\}|,\\
  L^*(N,G,\nu) \:\:\: = & \:\:\: |\{1\leq n,n'\leq N,n\neq n':\\
  &\:\:\: \exists j,j'\in\mathbb{Z}^d,1\leq ||j||_{\infty}|,||j'||_{\infty}\leq G,M_n^Tj\pm M_{n'}^Tj'=\nu\}|,\\
  L(N,G) \:\:\: = & \:\:\: \sup_{\nu\neq 0}L(N,G,\nu).
 \end{aligned}
\end{equation}
and for all $N\geq 1$, $G\geq 1$ and $\nu\in\mathbb{Z}$. Let $f:\mathbb{R}\to\mathbb{R}$ be some function of finite total variation which is one-periodic and satisfies $\mathbb{E}[f]=0$ as well as (\ref{2}) for some lacunary sequence satisfying the Hadamard gap condition (\ref{17}).
Aistleitner and Berkes showed that if for any fixed $G\geq 1$ we have $L(N,G)=o(N)$ for $N\to\infty$ then (\ref{1}) holds.

\subsection{Random Matrix ensembles}

Lacunary function systems have applications in random matrix theory. Random matrices are of particular interest in theoretic physics. Initially they were introduced by Wigner \cite{W55} to describe properties of atoms with heavy nuclei.
Consider an ensemble $X_N=(\frac{1}{\sqrt{N}}X_{n,n'})_{1\leq n,n'\leq N}$ of symmetric random $N\times N$ matrices such that except for the symmetry condition the entries $X_{n,n'}$ are independent random variables with mean zero, unit variance and universally bounded moments. Wigner showed that the mean empirical eigenvalue distribution converges weakly to the semicircle law as the size of the matrix tends to infinity. Being more precise in \cite{W58} he proved
\begin{equation}
 \label{252}
 \lim_{N\to\infty}\mathbb{E}\left[\frac{1}{N}\Tr\left(X_N^K\right)\right]=
 \begin{cases}
  C_{K/2}=\frac{K!}{(K/2)!(K/2+1)!}, & K \textnormal{ even},\\
  0, & K \textnormal{ odd}
 \end{cases}
\end{equation}
where the coefficients $C_K$ denote the Catalan numbers. Note that many problems in combinatorics have solutions which are related to Catalan numbers. For example, $C_{K/2}$ is the number of non-crossing pair partitions of $\{1,\ldots,K\}$ for some even $K$, i.e. there exist precisely $C_{K/2}$ pair partitions $\Delta$ such that there are no $k_1<k_2<k_3<k_4$ with $k_1,k_3\in\mathcal{S}$ and $k_2,k_4\in\mathcal{S}'$ for some $\mathcal{S},\mathcal{S}'\in\Delta$ with $\mathcal{S}\neq\mathcal{S}'$.
This fact plays an important part in the proof that in the limit the mean expectation the values of the traces of $X_N^K$ coincide with the moments of the semicircle law which has density $1/2\pi\cdot\sqrt{4-x^2}\cdot\mathbf{1}_{x^2\leq 4}$.\\ 
Having many other applications in physics, e.g. in quantum chaos or in telecommunications, in pure mathematics, e.g. in number theory, and further areas, random matrix models have been studied intensively in the last decades. In recent years the question arose whether the asymptotic behaviour still holds if the independence condition is weakened. Besides investigations on some specific models so far there are only some few attempts in this area. Schenker and Schulz-Baldes \cite{SS05} defined an ensemble of random matrices $X_N$ where for any $N$ there exists an equivalence relation $\sim_N$
on the set of entries of $X_N$ such that entries from different equivalence classes are independent while entries from the same class may be correlated. Observe that this model generalizes the classical case where the equivalence classes are of the form $\{X_{n,n'},X_{n',n}\}$. They showed that if the classes are not too large, i.e.
\begin{equation*}
 \max_{n\in\{1.\ldots,N\}}|\{n',m,m'\in\{1,\ldots,N\}:(n,n')\sim_N(m,m')\}|=o(N^2)
\end{equation*}
and
\begin{equation*}
 \max_{n,n',m\in\{1,\ldots,N\}}|\{m'\in\{1,\ldots,N\}:(n,n')\sim_N(m,m')\}|=C
\end{equation*}
for some absolute constant $C>0$, and if not too many different entries of the same class lie on the same row resp. the same column, i.e.
\begin{equation*}
 |\{n,n',m'\in\{1,\ldots,N\}:(n,n')\sim_N(n',m'),n\neq m'\}|=o(N^2),
\end{equation*}
then the mean empirical eigenvalue distribution converges weakly to the semicircle law. They proved the result without any further condition on the correlation of two entries in the same equivalence class. It is natural to consider random matrices such the correlation decays with the distance of two entries.
Therefore it is reasonable to study matrix models where each entry has except for some small errors only a finite range of dependence. As Anderson and Zeitouni \cite{AZ08} showed the converge to the semicircle law does not hold in general under this assumption, but further conditions are necessary. Thus although the
conditions given by Schenker and Schulz-Baldes appear to be too strict they can not be weakened without further constraints on the correlation of different entries. Friesen and L\"owe \cite{FL12} studied random matrix ensembles with stochastically independent diagonals but correlated entries on the diagonal. They showed convergence to the semicircle law
for
\begin{equation*}
 \max_{n,n'\in\{1,\ldots,N\}}\mathbb{E}[X_{n,n'}X_{n+t,n'+t}]\leq Ct^{-\varepsilon}
\end{equation*}
for some absolute constants $C>0$ and $\varepsilon>0$. Although it seems to be more reasonable to study ensembles with independent rows or columns rather than ensembles with independent diagonals these ensembles provide some difficulties. Since the symmetry condition is necessary to have real eigenvalues it implies that the columnwise independence and rowwise dependence turns into rowwise independence and
columnwise dependence by crossing the main diagonal. This not only seems to be not natural from a stochastic point of view, but also the mean empirical eigenvalue does not appear to converge to the semicircle law in general as simulations show.

\subsection*{Main result}




We discuss a random matrix ensemble where the entries are taken from a multivariate lacunary system and show that the mean spectral distribution of this ensemble converges weakly to the semicircle law.

\begin{theorem}
\label{253}
Let $(M_{n,1})_{n\geq 1}$ and $(M_{n,2})_{n\geq 0}$ be two integer-valued lacunary sequences satisfying the Hadamard gap condition (\ref{17}) such that there exist $C>0$ and $\varepsilon>0$ with
\begin{equation}
 \label{251}
 \sum_{\substack{n,n'\in\{1,\ldots,2N\},\\ n>n'}}\sum_{j,j'\in\{1,\ldots,N^K\}}(jj')^{-1}\mathbf{1}_{|jM_{n,1}-j'M_{n',1}|<1/2\cdot q^{-CN^{1-\varepsilon}}M_{n',1}}=o(N)
\end{equation}
for any $K\in\mathbb{N}$. Furthermore let $f:\mathbb{R}^2\to\mathbb{R}$ be some function of finite total variation in the sense of Hardy and Krause satisfying
\begin{equation}
 \label{221}
 f(x+z)=f(x) \textnormal{ for all }x\in [0,1)^2,z\in\mathbb{Z}^2, \quad \int_{[0,1)^2}f(x)\,dx=0, \quad \int_{[0,1)^2}f^2(x)\,dx=1.
\end{equation}
Additionally assume that $f$ satisfies
\begin{equation}
 \label{223}
 \int_{[0,1)}f(x_1,y_2)\,dx_1=\int_{[0,1)}f(y_1,x_2)\,dx_2=0
\end{equation}
for any fixed $y_1,y_2\in[0,1)$. Let the Fourier series of $f$ exist and converge to $f$.
Now define a symmetric random matrix ensemble $(X_N)_{N\geq 1}$ by setting $X_N=(X_{n,n'})_{1\leq n,n'\leq N}$ with
\begin{equation}
 \label{222}
 X_{n,n'}=\frac{1}{\sqrt{N}}f(M_{n+n',1}x_1,M_{|n-n'|,2}x_2)
\end{equation}
for any $n,n'\in\{1,\ldots,N\}$ and $N\in\mathbb{N}$. Then the mean empirical distribution of the eigenvalues converges almost surely weakly to the semicircle law, i.e.
\begin{equation}
 \label{250}
 \lim_{N\to\infty}\mathbb{E}\left[\frac{1}{N}\Tr\left(X_N^K\right)\right]=
 \begin{cases}
  \frac{K!}{(K/2)!(K/2+1)!}, & K \textnormal{ even},\\
  0, & K \textnormal{ odd}.
 \end{cases}
\end{equation}
\end{theorem}

\section{Preliminaries}

In this chapter we repeat some basic results on periodic functions of finite total variation resp. lacunary sequences which are going to be used in the subsequent chapters.

For some integer $d\geq 1$ set $I=\{1,\ldots,d\}$.
We now introduce the total variation in the sense of Hardy and Krause for periodic functions on $\mathbb{R}^d$.
Let $f:\mathbb{R}^d\to\mathbb{R}$ be some periodic function, i.e. $f$ satisfies $f(x+z)=f(x)$ for all $x\in\mathbb{R}^d$ and $z\in\mathbb{Z}^d$. For some subset $J\subseteq I$ and points $a,b\in [0,1)^{|J|}$ with $a_i\leq b_i$ for all $i\in J$ and some $z\in [0,1)^{|I\backslash J|}$ define
\begin{equation*}
 \Delta_J(f,a,b,z)=\sum_{\delta\in \{0,1\}^{|J|}}(-1)^{\sum_{i\in J}\delta_i}f(c_{\delta})
\end{equation*}
where $c_{\delta}=(c_{\delta,1},\ldots,c_{\delta,d})$ is defined by $c_{\delta,i}=\delta_i a_i+(1-\delta_i)b_i$ for $i\in J$ and $c_i=z_i$ for $i\notin J$. A finite set $\mathcal{Y}_i=\{y_1\ldots,y_{m(i)}\}\subset [0,1)$ with $0=y_1<\ldots<y_{m(i)}<1$ for some positive integer $m(i)$ is called a ladder. A multidimensional ladder on $[0,1)^d$ has the form $\mathcal{Y}=\prod_{i\in I}\mathcal{Y}_i$. For a multidimensional ladder $\mathcal{Y}$, a subset $J\subseteq I$ and $z\in[0,1)^d$ set $\mathcal{Y}_{J,z}=\prod_{i\in J}\mathcal{Y}_i\times \prod_{i\notin J}\{z_i\}$.
For $y\in\mathcal{Y}_{J,z}$ define $y_+\in [0,1)^{|J|}\times \prod_{i\notin J}\{z_i\}$ such that $y_{+,i}$ is the successor of $y_i$ in $\mathcal{Y}_i$ resp. $1$ if $y_i$ is the largest element in $\mathcal{Y}_i$. Then we define the variation of $f$ over $\mathcal{Y}_J$ by
\begin{equation*}
 V_{\mathcal{Y}_{J,z}}(f)=\sum_{y\in\mathcal{Y}_{J,z}}|\Delta_J(f,y,y_+,z)|.
\end{equation*}
Denote the set of all ladders $\mathcal{Y}_{J,z}$ by $\mathbb{Y}_{J,z}$. Then the total variation of $f$ over $[0,1)^{|J|}$ is defined by
\begin{equation*}
 V_{J}(f)=\sup_{z\in[0,1)^{|I|}}\sup_{\mathcal{Y}_{J,z}\in\mathbb{Y}_{J,z}}V_{\mathcal{Y}_{J,z}}(f).
\end{equation*}
The total variation of $f$ on $[0,1)^d$ in the sense of Hardy and Krause is
\begin{equation*}
 V_{HK}(f)=\sum_{J\subseteq I,J\neq \emptyset}V_{J}(f).
\end{equation*}
A function $f$ is called to be of finite total variation if $V_{HK}(f)<\infty$.

The following Lemma was proved in \cite{Z68}:

\begin{lemma}
\label{40}
 Let $f(x)=\sum_{j\in\mathbb{Z}^d\backslash\{0\}}a_j\cos(2\pi\langle j,x\rangle)+b_j\sin(2\pi\langle j,x\rangle)$ be a periodic function of finite total variation in the sense of Hardy and Krause. Then we have
\begin{equation*}
|a_j|,|b_j|\leq C\left(\prod_{i\in I,j_i\neq 0}\frac{1}{2\pi |j_i|}\right)V_{\{i\in I: j_i\neq 0\}}(f).
\end{equation*}
for some absolute constant $C>0$ and all $j\in\mathbb{Z}^d\backslash\{0\}$.
\end{lemma}

Observe that hereafter we always write $C$ for some absolute constant which may vary from line to line. Furthermore we always assume that $f\in L^2(\mathbb{R}^d,\mathbb{R})$ is a periodic function of finite total variation in the sense of Hardy and Krause.

For $\Gamma\in\mathbb{N}_0^d$ we denote the $\Gamma$th Dirichlet kernel by
\begin{equation}
\label{301}
D_{\Gamma}(x)=\sum_{j\in\mathbb{Z}^d,|j_i|\leq \Gamma_i}\cos(2\pi\langle j,x\rangle).
\end{equation}
Then the $\Gamma$th partial sum of $f$ is defined by
\begin{equation}
 \label{127}
 \psi_{\Gamma}(x)=\int_{[0,1)^d}f(x+t)D_{\Gamma}(t)\,dt=\sum_{j\in\mathbb{Z}^d,|j_i|\leq\Gamma_i}a_j\cos(2\pi\langle j,x\rangle)+b_j\sin(2\pi\langle j,x\rangle)
\end{equation}
for suitable numbers $a_j,b_j\in\mathbb{R}$ for all $j\in\mathbb{Z}^d$ with $|j_i|\leq \Gamma$ for all $i\in\{1,\ldots,d\}$. Set $\rho_{\Gamma}(x)=f(x)-\psi(x)$. If $\Gamma_i=G$ for all $i\in I$ we simply write $\psi_G(x)$ resp. $\rho_G(x)$. The $G$th Fej\'er mean of $f$ is defined by
\begin{equation}
 \label{92}
 \begin{aligned}
  p_G(x) \:\:\: & = \:\:\: \frac{1}{(G+1)^d}\sum_{\Gamma\in\mathbb{N}_0^d,||\Gamma||_{\infty}\leq G}\psi_{\Gamma}(x)\\
  & = \:\:\: \sum_{j\in\mathbb{Z}^d,|j_i|\leq G}a'_j\cos(2\pi\langle j,x\rangle)+b'_j\sin(2\pi\langle j,x\rangle)
 \end{aligned}
\end{equation}
where
\begin{equation*}
 a'_j=a_j\prod_{i\in I}\frac{G+1-|j_i|}{G+1}, \quad  b'_j=b_j\prod_{i\in I}\frac{G+1-|j_i|}{G+1}
\end{equation*}
for all $j\in\mathbb{Z}^d$ with $||j||_{\infty}\leq G$. Observe that
\begin{equation*}
 p_G(x)=\int_{[0,1)^d}f(x+t)K_G(t)\,dt
\end{equation*}
where $K_G(t)=K_{\Gamma}(t)$ with $\Gamma_i=G$ for all $i\in I$ and $K_{\Gamma}(t)=\prod_{i\in I}K_{\Gamma_i}(t_i)=\prod_{i\in I}K_{G}(t_i)$ is the $d$-dimensional $G$th Fej\'er kernel and $K_G(t_i)$ is the one-dimensional $G$th Fej\'er kernel defined by
\begin{eqnarray*}
 K_G(t_i) & = & \frac{1}{G+1}\sum_{l=0}^{G}\sum_{j_i=-l}^{l}\cos(2\pi\langle j_i,x\rangle)\\
 & = & \frac{1}{G+1}\sum_{l=0}^{G}\frac{\sin(2\pi\langle l+1/2,t_i \rangle)}{\sin(2\pi\langle 1/2,t_i \rangle)}\\
 & = & \frac{1}{G+1}\frac{(\sin(2\pi\langle (G+1)/2,t_i \rangle))^2}{2(\sin(2\pi\langle 1/2,t_i \rangle))^2}\\
 & \geq & 0.
\end{eqnarray*}
Therefore we have
\begin{eqnarray*}
 |p_G(x)| & = & \left|\int_{[0,1)^d}f(x+t)K_G(t)\,dt\right|\\
 & \leq & ||f||_{\infty}\left|\int_{[0,1)^d}K_G(t)\,dt\right|\\
 & \leq & ||f||_{\infty}.
\end{eqnarray*}
Now we define $r_G(x)=f(x)-p_G(x)$.

\begin{lemma}
 \label{94}
 Let $f(x)=\sum_{j\in\mathbb{Z}^d\backslash\{0\}}a_j\cos(2\pi\langle j,x\rangle)+b_j\sin(2\pi\langle j,x\rangle)$ be some periodic function satisfying $V_{HK}(f)\leq 1$. Then there exists some absolute constant $C>0$ such that for any $G\geq d$ the function $r_G$ satisfies
 \begin{equation*}
  ||r_G||_2^2\leq CdG^{-1}.
 \end{equation*}
\end{lemma}

\textit{Proof.} 
Observe that for any $\varepsilon=\varepsilon(d,G)>0$ there exists some trigonometric polynomial $f'$ such that $||f-f'||_2\leq \varepsilon$. Let $p_G'$ be the $G$th Fej\'er mean of $f'$. We obtain $||p_G-p_G'||_2\leq \varepsilon$. Therefore we have $||r_G||_2\leq ||f'-p_g'||_2+C\sqrt{dG^{-1}}$. Thus it is enough to prove the statement of the Lemma for trigonometric polynomials $f$.
Set
\begin{equation}
  \label{83}
 r_G(x)=\sum_{j\in\mathbb{Z}^d\backslash\{0\}}\tilde{a}_j\cos(2\pi\langle j,x\rangle)+\tilde{b}_j\sin(2\pi\langle j,x\rangle)
\end{equation}
where
\begin{equation*}
 \tilde{a}_j=
 \begin{cases}
 a_j, & ||j||_{\infty}>G,\\
 a_j\left(1-\prod_{i\in I}\left(1-\frac{|j_i|}{G+1}\right)\right), & ||j||_{\infty}\leq G
 \end{cases}
\end{equation*}
and $\tilde{b}_j$ is defined analogously.
We have
\begin{equation}
  \label{91}
 ||r_G||_2^2=\sum_{j\in\mathbb{Z}^d\backslash\{0\}}\tilde{a}_j^2+\tilde{b}_j^2=\underbrace{\sum_{0<||j||_{\infty}\leq G}\tilde{a}_j^2+\tilde{b}_j^2}_{(*)}+\underbrace{\sum_{||j||_{\infty}>G}\tilde{a}_j^2+\tilde{b}_j^2}_{(**)}.
\end{equation}
For some given nonempty $J\subseteq I$ set 
\begin{eqnarray*}
 D(G,J) & = & \{j\in\mathbb{Z}^d:1\leq |j_i|\leq G\textnormal{ for } i\in J,j_i=0\textnormal{ for } i\notin J\},\\
 D'(G,J) & = & \{j\in\mathbb{Z}^d:j_i\neq 0 \textnormal{ for }i\in J,j_i=0\textnormal{ for } i\notin J,j\notin D(G,J)\}.
\end{eqnarray*}
To estimate $(*)$ we first by Lemma \ref{40} observe 
\begin{equation*}
 \sum_{0<||j||_{\infty}\leq G}\tilde{a}_j^2+\tilde{b}_j^2 \leq  2\sum_{J\subseteq I,J\neq \emptyset}\sum_{j\in D(G,J)}\left(\prod_{i\in J}\frac{1}{2\pi|j_i|}\right)^2\left(1-\prod_{i\in J}\left(1-\frac{|j_i|}{G+1}\right)\right)^2V_J(f)^2.
\end{equation*}
By definition of $V_{HK}(f)$ it is enough to show
\begin{equation}
 \label{89}
 V(G,J)=2\sum_{j\in D(G,J)}\left(\prod_{i\in J}\frac{1}{2\pi|j_i|}\right)^2\left(1-\prod_{i\in J}\left(1-\frac{|j_i|}{G+1}\right)\right)^2\leq CG^{-1}
\end{equation}
for some absolute constant $C>0$. By decomposing we have
\begin{equation*}
 \begin{aligned}
 V(G,J)=2\sum_{K,K'\subseteq J,K,K'\neq 0}\sum_{j\in D(G,J)} & \prod_{i\in K}\frac{1}{2\pi|j_i|}\frac{|j_i|}{G+1}\prod_{i\in J\backslash K}\frac{1}{2\pi|j_i|}\left(1-\frac{|j_i|}{G+1}\right)\\
 \cdot & \prod_{i\in K'}\frac{1}{2\pi|j_i|}\frac{|j_i|}{G+1}\prod_{i\in J\backslash K'}\frac{1}{2\pi|j_i|}\left(1-\frac{|j_i|}{G+1}\right).
 \end{aligned}
\end{equation*}
Thus we get
\begin{equation*}
 V(G,J)=2\sum_{K,K'\subseteq J,K,K'\neq 0}W_1(K,K')\cdot W_2(K,K')\cdot W_3(K,K')\cdot W_4(K,K')
\end{equation*}
where
\begin{eqnarray*}
 W_1(K,K') & = & \prod_{i\in K\cap K'}\frac{1}{(2\pi)^2}\sum_{j_i=-G}^G\frac{1}{(G+1)^2}\leq \prod_{i\in K\cap K'}\frac{1}{4G},\\
 W_2(K,K') & = & \prod_{i\in K\cap J\backslash K'}\frac{1}{(2\pi)^2}\sum_{j_i=-G}^G\frac{1}{G+1}\left(\frac{1}{|j_i|}-\frac{1}{G+1}\right)\leq \prod_{i\in K\cap J\backslash K'}\frac{\log(G)}{4G},\\
 W_3(K,K') & = & \prod_{i\in J\backslash K\cap K'}\frac{1}{(2\pi)^2}\sum_{j_i=-G}^G\left(\frac{1}{|j_i|}-\frac{1}{G+1}\right)\frac{1}{G+1}\leq \prod_{i\in J\backslash K\cap K'}\frac{\log(G)}{4G},\\
 W_4(K,K') & = & \prod_{i\in J\backslash K\cap J\backslash K'}\frac{1}{(2\pi)^2}\sum_{j_i=-G}^G\left(\frac{1}{|j_i|}-\frac{1}{G+1}\right)^2\leq \prod_{i\in J\backslash K\cap J\backslash K'}\frac{1}{4}.
\end{eqnarray*}
Since $K\cap K'\neq\emptyset$ or $K\cap J\backslash K'\neq\emptyset$ and $J\backslash K\cap K'\neq\emptyset$ we conclude
\begin{equation}
 \label{90}
 V(G,J)\leq \sum_{K,K'\subseteq J,K,K'\neq 0}\frac{C}{4^{|J|}G}\leq CG^{-1}
\end{equation}
for some absolute constant $C>0$ and therefore (\ref{89}) is verified.
We now estimate $(**)$. By Lemma \ref{40} we have
\begin{eqnarray*}
 \sum_{||j||_{\infty}>G}\tilde{a}_j^2+\tilde{b}_j^2 & \leq & \sum_{J\subseteq I,J\neq\emptyset}\sum_{j\in D'(G,J)}a_j^2+b_j^2\\
 & \leq & 2\sum_{J\subseteq I,J\neq\emptyset}\sum_{j\in D'(G,J)}\left(\prod_{i\in J}\frac{1}{(2\pi|j_i|)^2}\right)V_J(f)^2.
\end{eqnarray*}
We furthermore for some nonempty $J\subseteq I$ get
\begin{eqnarray*}
 \sum_{j\in D'(G,J)}\left(\prod_{i\in J}\frac{1}{(2\pi|j_i|)^2}\right) & \leq &  \sum_{l=1}^{|J|}{|J|\choose l}\left(2\sum_{j=1}^G\frac{1}{(2\pi j)^2}\right)^{|J|-l}\cdot \left(2\sum_{j=G+1}^{\infty}\frac{1}{(2\pi j)^2}\right)^{l}\\
 & \leq & \sum_{l=1}^{|J|}{|J|\choose l}\left(\frac{1}{2\pi^2G}\right)^l\\
 & \leq & \sum_{l=1}^{|J|}\left(\frac{d}{2\pi^2G}\right)^l\\
 & \leq & CdG^{-1}
\end{eqnarray*}
for some absolute constant $C>0$. Therefore we have
\begin{equation*}
 \sum_{||j||_{\infty}>G}\tilde{a}_j^2+\tilde{b}_j^2\leq CdG^{-1}
\end{equation*}
and the Lemma is proved.

\vspace{3ex}

With $L(N,G,\nu)$ as defined in (\ref{22}) we have

\begin{lemma}
\label{102}
 Let $f:\mathbb{R}^d\to\mathbb{R}$ be a periodic function satisfying $V_{HK}(f)\leq 1$ and let $(M_n)_{n\geq 0}$ be a lacunary sequence of matrices satisfying the Hadamard gap condition (\ref{17}). Then we have
 \begin{equation}
  \label{15}
  \mathlarger{\mathlarger{\int}}_{[0,1)^d}\left(\sum_{n=1}^Nf(M_nx)\right)^2\,dx\leq C(\log(d)||f||_2^2+||f||_2)N
 \end{equation}
 where $C>0$ is an absolute constant depending only on $q$. If the sequence $(M_n)_{n\geq 1}$ furthermore satisfies $L(N,G,0)=o(N)$
 for any fixed $G\geq 1$ then we have
 \begin{equation}
 \label{96}
 \lim_{N\to\infty}\frac{1}{N}\mathlarger{\mathlarger{\int}}_{[0,1)^d}\left(\sum_{n=1}^Nf(M_nx)\right)^2\,dx=||f||_2^2.
 \end{equation}
\end{lemma}

\vspace{3ex}

Note that hereafter we write $\log(x)$ for $\max(1,\log(x))$.

\vspace{3ex}

\textit{Proof.} For $||j||_{\infty},||j'||_{\infty}\leq G$ and $k>\log_q(G)$ we have
\begin{equation}
  \label{93}
 ||M^T_nj'||_{\infty}\leq ||M^T_n||_{\infty}||j'||_{\infty}\leq q^{-k}||M^T_{n+k}j||_{\infty}||j'||_{\infty}< ||M^T_{n+k}j||_{\infty}.
\end{equation}
Therefore we obtain $M^T_nj'\neq M^T_{n+k}j$. Now let $p_G$ be the $G$th Fej\'er mean of $f$. Then for some $k>\log_q(G)$ we have
\begin{equation*}
 p_G(M_nx)p_G(M_{n+k}x)=\sum_u\alpha_u\cos(2\pi\langle u,x\rangle)+\beta_u\sin(2\pi\langle u,x\rangle)
\end{equation*}
where any $u$ is of the form $M^T_nj\pm M^T_{n+k}j'$ for some $1\leq ||j||_{\infty},||j'||_{\infty}\leq G$. Therefore by (\ref{93}) we get
\begin{equation}
  \label{95}
 \int_{[0,1)^d}p_G(M_nx)p_G(M_{n+k}x)\,dx=0
\end{equation}
for $k>\log_q(G)$. By Lemma \ref{94} for any $k\geq 1$ there is a trigonometric polynomial $g_k$ with $dq^{2k}-1<\deg(g_k)\leq dq^{2k}$ such that
\begin{equation*}
 ||f-g_k||_2\leq Cq^{-k}.
\end{equation*}
Therefore for $k'>\log_q(dq^{2k})$ by (\ref{95}) and Cauchy-Schwarz inequality we have
\begin{equation}
\label{18}
\begin{aligned}
 \left|\int_{[0,1)^d} f(M_{n}x)f(M_{n+k'}x)\,dx\right| \:\:\: \leq & \:\:\: \left|\int_{[0,1)^d} (f-g_k)(M_{n}x)f(M_{n+k'}x)\,dx\right|\\
 &\:\:\: +\left|\int_{[0,1)^d} g_k(M_{n}x)g_k(M_{n+k'}x)\,dx\right|\\
 &\:\:\: +\left|\int_{[0,1)^d} g_k(M_{n}x)(f-g_k)(M_{n+k'}x)\,dx\right|\\
 \leq &\:\:\: 2C||f||_2q^{-k}
 \end{aligned}
\end{equation}
since $||g_k||_2\leq ||f||_2$. We obtain
\begin{multline*}
\mathlarger{\mathlarger{\int}}_{[0,1)^d}\left(\sum_{n=1}^Nf(M_nx)\right)^2\,dx\\
\begin{aligned}
\leq & N||f||_2^2+2\sum_{n=1}^N\sum_{k'=1}^{N-n}\left|\int_{[0,1)^d}f(M_nx)f(M_{n+k'}x)\,dx\right|\\
\leq & N||f||_2^2+2N(\log_q(d)+2)||f||_2^2+2\sum_{n=1}^N\sum_{k'=\log_q(d)+3}^{N-n}\left|\int_{[0,1)^d}f(M_nx)f(M_{n+k'}x)\,dx\right|\\
\leq & N||f||_2^2+CN\log_q(d)||f||_2^2+\sum_{n=1}^N\sum_{k=1}^{\infty}C||f||_2q^{-k}\\
\leq & C(\log(d)||f||_2^2+||f||_2)N.
\end{aligned}
\end{multline*}
Thus (\ref{15}) is shown.\\
The proof of (\ref{96}) is similar. For $k'>\log_q(q^{2k})$ instead of (\ref{18}) we get
\begin{equation*}
 \left|\int_{[0,1)^d} f(M_{n}x)f(M_{n+k}x)\,dx\right| \leq C||f||_2d^{1/2}q^{-k}
\end{equation*}
and similarly
\begin{equation*}
 \left|\int_{[0,1)^d} s_1(M_{n}x)s_2(M_{n+k}x)\,dx\right| \leq C||f||_2d^{1/2}q^{-k}
\end{equation*}
where for $i\in\{1,2\}$ the function $s_i$ is of the form $p_{G_i}$ or $r_{G_i}$ for some suitable number $G_i>0$.\\
For any $G\geq 1$ we obtain
\begin{multline*}
 \left|\frac{1}{N}\mathlarger{\mathlarger{\int}}_{[0,1)^d}\left(\sum_{n=1}^Nf(M_nx)\right)^2\,dx-||f||_2^2\right|\\
 \begin{aligned}
 \leq & \frac{2}{N}\sum_{n=1}^N\sum_{k=1}^{N-n}\left|\int_{[0,1)^d}f(M_nx)f(M_{n+k}x)\,dx\right|\\
 \leq & \frac{C}{N}\sum_{n=1}^N\sum_{k=1}^{N-n}||f||_2d^{1/2}\min(G^{-1/2},q^{-k})\\
 & +\frac{2}{N}\sum_{n=1}^N\sum_{k=1}^{N-n}\left|\int_{[0,1)^d}p_G(M_nx)p_G(M_{n+k}x)\,dx\right|\\
 \leq & C||f||_2d^{1/2}G^{-1/2}\\
 & + (2G+1)^{2d}\frac{h(N)}{N}
 \end{aligned}
\end{multline*}
for some function $h(N)$ with $h(N)/N\to 0$. Observe that such a function exists by assumption on $L(N,G,0)$. Since the constant $G\geq 1$ can be chosen arbitrary, (\ref{96}) is shown.

\vspace{5ex}

\section{Proof of Theorem \ref{253}}
\label{S51}

For $K,N\in\mathbb{N}$ let $p=p_{N^K}$ be the $N^K$th Fej\'er mean of $f$ as defined in (\ref{92}) and set $r=r_{N^K}=f-p$. By Lemma \ref{94} we have $||r||_2^2\leq CN^{-K}$ for some absolute constant $C>0$.
We define the random matrix $(\tilde{X}_{n,n'})_{1\leq n,n'\leq N}$ by
\begin{equation}
 \label{226}
 \tilde{X}_{n,n'}=\frac{1}{\sqrt{N}}p(M_{n+n',1}x_1,M_{|n-n'|,2}x_2)
\end{equation}
for any $n,n'\in\{1,\ldots,N\}$.
Now we claim
\begin{equation}
 \label{225}
 \lim_{N\to\infty}\mathbb{E}\left[\frac{1}{N}\Tr\left(X_N^K\right)\right]=\lim_{N\to\infty}\mathbb{E}\left[\frac{1}{N}\Tr\left(\tilde{X}_N^K\right)\right].
\end{equation}
It is easy to see that
\begin{equation*}
 \begin{aligned}
  \mathbb{E}\left[\frac{1}{N}\Tr\left(X_N^K\right)\right] \:\:\: & = \:\:\: \mathbb{E}\left[\frac{1}{N}\sum_{n_1,\ldots,n_K=1}^N\prod_{k=1}^KX_{n_{k},n_{k+1}}\right]\\
  & = \:\:\: \mathbb{E}\left[\frac{1}{N^{K/2+1}}\sum_{n_1,\ldots,n_K=1}^N\prod_{k=1}^Kf(M_{n_k+n_{k+1},1}x_1,M_{|n_k-n_{k+1}|,2}x_2)\right]
 \end{aligned}
\end{equation*}
where $n_{1}=n_{K+1}$.
By decomposing $f=p+r$ we observe
\begin{multline}
 \label{227}
 \mathbb{E}\left[\frac{1}{N}\Tr\left(X_N^K\right)\right]\\
 \begin{aligned}
  = & \frac{1}{N^{K/2+1}}\sum_{n_1,\ldots,n_K=1}^N\mathbb{E}\left[\prod_{k=1}^Kp(M_{n_k+n_{k+1},1}x_1,M_{|n_k-n_{k+1}|,2}x_2)\right]\\
  & +\frac{1}{N^{K/2+1}}\sum_{n_1,\ldots,n_K=1}^N\sum_{\substack{J\subseteq\{1,\ldots,K\},\\ J\neq\emptyset}}\mathbb{E}\left[\prod_{k\in J}r(M_{n_k+n_{k+1},1}x_1,M_{|n_k-n_{k+1}|,2}x_2)\right.\\
  & \quad\quad\quad\quad\quad\quad\quad\quad\quad\quad\quad\quad\quad\quad\,\cdot\left.\prod_{k\notin J}p(M_{n_k+n_{k+1},1}x_1,M_{|n_k-n_{k+1}|,2}x_2)\right].
 \end{aligned}
\end{multline}
Since any set $J$ is nonempty there exists some $k_J\in J$ for each $J$. Thus by Cauchy-Schwarz inequality we get
\begin{multline*}
 \left|\mathbb{E}\left[\prod_{k\in J}r(M_{n_k+n_{k+1},1}x_1,M_{|n_k-n_{k+1}|,2}x_2)\prod_{k\notin J}p(M_{n_k+n_{k+1},1}x_1,M_{|n_k-n_{k+1}|,2}x_2)\right]\right|\\
 \begin{aligned}
  \leq & \mathbb{E}\left[r^2(M_{n_{k_J}+n_{k_J+1},1}x_1,M_{|n_{k_J}-n_{k_J+1}|,2}x_2)\right]^{1/2}\\
  & \cdot\mathbb{E}\left[\prod_{k\in J\backslash\{k_J\}}r^2(M_{n_k+n_{k+1},1}x_1,M_{|n_k-n_{k+1}|,2}x_2)\prod_{k\notin J}p^2(M_{n_k+n_{k+1},1}x_1,M_{|n_k-n_{k+1}|,2}x_2)\right]^{1/2}
 \end{aligned}
\end{multline*}
and therefore we obtain
\begin{multline*}
 \left|\mathbb{E}\left[\prod_{k\in J}r(M_{n_k+n_{k+1},1}x_1,M_{|n_k-n_{k+1}|,2}x_2)\prod_{k\notin J}p(M_{n_k+n_{k+1},1}x_1,M_{|n_k-n_{k+1}|,2}x_2)\right]\right|\\
 \begin{aligned}
  \leq & CN^{-K/2}\left(\prod_{k\in J\backslash\{k_J\}}||r||_{\infty}^2\right)^{1/2}\left(\prod_{k\notin J}||p||_{\infty}^2\right)^{1/2}\\
  \leq & CN^{-K/2}||f||_{\infty}^{K-1}
 \end{aligned}
\end{multline*}

for some constant $C>0$ which only depends on $K$. Plugging this into (\ref{227}) we have
%
\begin{equation*}
 \begin{aligned}
  \mathbb{E}\left[\frac{1}{N}\Tr\left(X_N^K\right)\right] \:\:\: = & \:\:\: \frac{1}{N^{K/2+1}}\sum_{n_1,\ldots,n_K=1}^N\mathbb{E}\left[\prod_{k=1}^Kp(M_{n_k+n_{k+1},1}x_1,M_{|n_k-n_{k+1}|,2}x_2)\right]\\
  & \:\:\: + \frac{1}{N^{K/2+1}}\sum_{n_1,\ldots,n_K=1}^N\sum_{\substack{J\subseteq\{1,\ldots,K\},\\ J\neq \emptyset}}CN^{-K/2}||f||_{\infty}^{K-1}.
 \end{aligned}
\end{equation*}
Hence (\ref{225}) follows immediately. Thus it is enough to study the asymptotic behaviour of $\mathbb{E}[1/N\cdot\Tr(\tilde{X}_N^K)]$.
We have
\begin{equation}
 \label{224}
 f(x)=\sum_{j\in(\mathbb{Z}\backslash \{0\})^2}a_j\cos(2\pi\langle j,x\rangle)+b_j\sin(2\pi\langle j,x\rangle)
\end{equation}
for suitable numbers $a_j,b_j$ where without loss of generality we may assume $b_j=0$ for all $j\in(\mathbb{Z}\backslash\{0\})^2$. The general case is similar. Hence by definition we obtain
\begin{multline}
 \label{247}
 \mathbb{E}\left[\frac{1}{N}\Tr\left(\tilde{X}_N^K\right)\right]\\
 \begin{aligned}
  = & \frac{1}{N^{K/2+1}}\sum_{n_1,\ldots,n_K\in\{1,\ldots,N\}}\mathbb{E}\left[\prod_{k=1}^Kp(M_{n_k+n_{k+1},1}x_1,M_{|n_k-n_{k+1}|,2}x_2)\right]\\
  = & \sum_{\delta_1,\ldots,\delta_K\in\{0,1\}}\sum_{\substack{j_1,\ldots,j_K\in(\mathbb{Z}\backslash\{0\})^2,\\ ||j_k||_{\infty}\leq N^K\forall k}}\sum_{n_1,\ldots,n_K\in\{1,\ldots,N\}}\frac{1}{2^K}\frac{1}{N^{K/2+1}}\prod_{k=1}^Ka'_{j_k}\\
  & \cdot\mathbb{E}\left[\cos\left(2\pi\left(\sum_{k=1}^K(-1)^{\delta_k}j_{k,1}M_{n_k+n_{k+1},1}x_1+\sum_{k=1}^K(-1)^{\delta_k}j_{k,2}M_{|n_k-n_{k+1}|,2}x_2\right)\right)\right].
 \end{aligned}
\end{multline}
Now we decompose the system of solutions for $\sum_{k=1}^K(-1)^{\delta_k}j_{k,1}M_{n_k+n_{k+1},1}=0$ and $\sum_{k=1}^K(-1)^{\delta_k}j_{k,2}M_{|n_k-n_{k+1}|,2}=0$ into different sets. Therefore we rearrange the sequence $(M_{n_k+n_{k+1},1})_{k\in\{1,\ldots,K\}}$ resp. $(M_{|n_k-n_{k+1}|,2})_{k\in\{1,\ldots,K\}}$ in decreasing order.
Let $\pi_1$ and $\pi_2$ be permutations of $\{1,\ldots,K\}$ such that we have $n_{\pi_1(k)}+n_{\pi_1(k)+1}\geq n_{\pi_1(k')}+n_{\pi_1(k')+1}$ resp. $|n_{\pi_2(k)}-n_{\pi_2(k)+1}|\geq |n_{\pi_2(k')}-n_{\pi_2(k')+1}|$ for any $k\leq k'$. For $i\in\{1,2\}$ we define sequences $(l_{k,i})_{1\in\{1,\ldots,K\}}$ by
setting $l_{k,1}=n_{\pi_1(k)}+n_{\pi_1(k)+1}$ and $l_{k,2}=|n_{\pi_2(k)}-n_{\pi_2(k)+1}|$ for any $k\in\{1,\ldots,K\}$. For any two sets $J_1,J_2\subset\{1,\ldots,K\}$ let the set $\mathcal{A}_{J_1,J_2}$ consist of any solution
$\delta_1,\ldots,\delta_K,j_1,\ldots,j_K,n_1,\ldots,n_K$ such that for any $i\in\{1,2\}$ we have
\begin{equation}
 \label{245}
 \left|\sum_{k=1}^H(-1)^{\delta_{\pi_i(k)}}j_{\pi_i(k),i}M_{l_{k,i},i}\right|<\frac{1}{2}M_{l_{H,i},i}
\end{equation}
if and only if $H\in J_i$.
If $H\notin J_i$ for some $H\in\{1,\ldots,K\}$ then we have
\begin{equation*}
 \left|\sum_{k=H+1}^K(-1)^{\delta_{\pi_i(k)}}j_{\pi_i(k),i}M_{l_{k,i},i}\right|=\left|\sum_{k=1}^H(-1)^{\delta_{\pi_i(k)}}j_{\pi_i(k),i}M_{l_{k,i},i}\right|\geq\frac{1}{2}M_{l_{H,i},i}.
\end{equation*}
Thus we get
\begin{equation*}
 KN^KM_{l_{H+1,i},i}\geq \frac{1}{2}M_{l_{H,i},i}.
\end{equation*}
Simple calculation shows
\begin{equation}
 \label{229}
 l_{H,i}-l_{H+1,i}\leq \log_q(2K)+K\log_q(N).
\end{equation}
Now let $H\in J_i$ for some $H\geq 2$. We have
\begin{eqnarray*}
 \left|\sum_{k=1}^{H-1}(-1)^{\delta_{\pi_i(k)}}j_{\pi_i(k),i}M_{l_{k,i},i}\right| & \geq & \left|j_{\pi_i(H),i}M_{l_{H,i},i}\right|-\left|\sum_{k=1}^{H}(-1)^{\delta_{\pi_i(k)}}j_{\pi_i(k),i}M_{l_{k,i},i}\right|\\
 & > & \frac{1}{2}M_{l_{H,i},i}.
\end{eqnarray*}
Assume that there exists another solution $\delta'_1,\ldots,\delta'_K,j'_{1,i},\ldots,j'_{K,i},n'_1,\ldots,n'_K$ in $\mathcal{A}_{J_1,J_2}$ with $\delta'_{\pi_i(k)}=\delta_{\pi_i(k)}$,$j'_{\pi_i(k),i}=j_{\pi_i(k),i}$ and $l'_{k,i}=l_{k,i}$ for $k\in\{1,\ldots,H-1\}$ where any $l'_{k,i}$ is defined analogously to $l_{k,i}$.
Without loss of generality we may assume that $l_{H,i}\geq l'_{H,i}$. Then similar calculation as above shows
\begin{equation}
 \label{230}
 l_{H,i}-l'_{H,i}<\log_q(2K)+K\log_q(N).
\end{equation}
Now we are going to estimate
\begin{equation}
 \label{228}
 \sum_{\substack{\delta_1,\ldots,\delta_K,j_1,\ldots,j_K,\\ n_1,\ldots,n_K\in\mathcal{A}_{J_1,J_2}}}\frac{1}{2^K}\frac{1}{N^{K/2+1}}\prod_{k=1}^Ka'_{j_k}
\end{equation}
for any particular pair of sets $J_1,J_2$. Therefore at first we encounter all possible choices for $l_{1,i},\ldots,l_{K,i}$ and $i\in\{1,2\}$. Observe that by (\ref{229}) and (\ref{230}) the number of choices for $l_{k,i}$ is bounded by $C\log(N)$ for some constant $C>0$ independent of $N$ if $k-1\notin J_i$ or $k\in J_i$, otherwise it is bounded by $2N$. Now we assume that $J_1$ or $J_2$ is not $\{2,4,\ldots,K\}$ for an even integer $K$. Then the total number of choices for $l_{1,i},\ldots,l_{K,i}$ is bounded by $CN^{\lfloor K/2\rfloor -1}\log(N)^{\lceil K/2\rceil +1}$ for some constant $C>0$ independent of $N$. Furthermore there are at most $N$ choices for $n_1$ and together with $l_{1,i},\ldots,l_{K,i}$ this uniquely determines any $n_k$ for $k\in\{2,\ldots,K\}$. Thus the total number of choices for $n_1,\ldots,n_K$ is bounded by $CN^{\lfloor K/2\rfloor}\log(N)^{\lceil K/2\rceil +1}$. 
By Lemma \ref{40} we have
\begin{equation}
 \label{232}
 \left|\sum_{\substack{\delta_1,\ldots,\delta_K,\\ j_1,\ldots,j_K}}\frac{1}{2^K}\prod_{k=1}^Ka'_{j_k}\right|\leq C\log(N)^{K}
\end{equation}
for some constant $C>0$ independent of $N$. We conclude
\begin{equation}
 \label{239}
 \left|\sum_{\substack{(\delta_1,\ldots,\delta_K,j_1,\ldots,j_K,\\ n_1,\ldots,n_K)\in\mathcal{A}_{J_1,J_2}}}\frac{1}{2^K}\frac{1}{N^{K/2+1}}\prod_{k=1}^Ka'_{j_k}\right|=o(1)
\end{equation}
if $J_1$ or $J_2$ is not $\{2,4,\ldots,K\}$ for an even integer $K$. Thus from now on we may assume 
\begin{equation}
 \label{231}
 J_1=J_2=\{2,4,\ldots,K\}
\end{equation}
and we simply write $\mathcal{A}_{\{2,4,\ldots,K\}}=\mathcal{A}_{J_1,J_2}$.
We further may assume that for  any $C>0$ and $\varepsilon>0$ we have 
\begin{equation}
 \label{238}
 l_{k,i}\geq l_{k+1,i}+2CN^{1-\varepsilon}
\end{equation}
for any $k\in\{2,4,\ldots,K-2\}$ and sufficiently large $N$ where $C$ denotes the constant used in (\ref{251}) since repeating the above argumentation reveals that all other solutions may be neglected.
Observe that $\pi_1$ and $\pi_2$ define decompositions $\Delta_{\pi_1}$ and $\Delta_{\pi_2}$ of $\{1,\ldots,K\}$ into pairs such that for $i\in\{1,2\}$ any pair $\{k,k'\}\in\Delta_{\pi_i}$ satisfies $\pi_i(k)+1=\pi_i(k')$ with $\pi_i(k)$ odd. We now claim that all solutions $\delta_1,\ldots,\delta_K$, $j_1,\ldots,j_K$, $n_1,\ldots,n_K$ such that $\Delta_{\pi_1}\neq \Delta_{\pi_2}$ may be neglected. Therefore we define a decomposition $\Delta$ of $\{1,\ldots,K\}$ such that for any $k,k'$ in different subsets we
have $\{k,k'\}\notin \Delta_{\pi_i}$ for any $i\in\{1,2\}$. Observe that for $\Delta_{\pi_1}\neq \Delta_{\pi_2}$ there are at most $K/2-1$ subsets. Now we encounter all solutions. The number of possibilities to decompose $\{1,\ldots,K\}$ into pairs in two different ways is independent of $N$. Thus we may assume that $\Delta$ is fixed. We determine the $n_k$ with $k\in\{1,\ldots,K\}$ in increasing order of the index $k$. There are $N$ choices for $n_1$. Furthermore any choice of $n_k$ and $n_{k+1}$ uniquely defines $l_{\pi_i^{-1}(k),i}$ for $i\in\{1,2\}$.  If $\{1,\ldots,k-1\}\cap\mathcal{S}\neq\emptyset$ with $k\in\mathcal{S}$ for $\mathcal{S}\in\Delta$
then by (\ref{229}) and (\ref{231}) the number of choices for $n_{k+1}$ is bounded by $C_1\log(N)^{C_2}$ for some constants $C_1,C_2>0$ independent of $N$. Otherwise the number of choices is bounded by $N$. Hence the total number of choices for $n_1,\ldots,n_K$ is bounded by $C_1N^{|\Delta|+1}\log(N)^{C_2}$ for some constants $C_1,C_2>0$ independent of $N$. Together with (\ref{232}) we obtain
\begin{equation}
 \label{240}
 \left|\sum_{\substack{(\delta_1,\ldots,\delta_K,j_1,\ldots,j_K,\\ n_1,\ldots,n_K)\in\mathcal{A}_{\{2,4,\ldots,K\}},\\ \Delta_{\pi_1}\neq \Delta_{\pi_2}}}\frac{1}{2^K}\frac{1}{N^{K/2+1}}\prod_{k=1}^Ka'_{j_k}\right|=o(1).
\end{equation}
Thus we restrict ourselves to the case $\Delta_{\pi_1}=\Delta_{\pi_2}$. 
Now we claim that pair partitions which are not non-crossing may be neglected. Thus we encounter all pair partitions with $\{k_1,k_3\},\{k_2,k_4\}\in\Delta$ for some $k_1<k_2<k_3<k_4$. We begin by counting the choices for $n_{k_4+1},\ldots,n_K$, $n_{K+1}=n_1,\ldots,n_{k_2}$. By (\ref{229}) the number of choices for $l_{k_3,1}$ and $l_{k_3,2}$ is bounded by $C\log(N)$ with some constant $C>0$ independent of $N$ for any fixed $n_{k_1}$ and $n_{k_1+1}$. Consequently, the number of choices for $n_{k_3}$ is bounded by $C\log(N)^2$ for some constants $C>0$ independent of $N$. Therefore we now determine $n_{k_3-1},\ldots,n_{k_2+1}$ in decreasing order and $n_{k_3+1},\ldots,n_{k_4}$ in increasing order.
Thus $l_{k_2}$ and $l_{k_4}$ already are uniquely defined by $n_1$ and $(l_{k,1})_{k\in\{1,\ldots,K\}\backslash\{k_2,k_4\}}$. Then using a similar argumentation as above we have
\begin{equation}
 \label{241}
 \left|\sum_{\substack{(\delta_1,\ldots,\delta_K,j_1,\ldots,j_K,\\ n_1,\ldots,n_K)\in\mathcal{A}_{\{2,4,\ldots,K\}},\\ \Delta \notin\mathcal{D}}}\frac{1}{2^K}\frac{1}{N^{K/2+1}}\prod_{k=1}^Ka'_{j_k}\right|=o(1)
\end{equation}
where $\mathcal{D}$ denotes the set of all non-crossing pair partitions of $\{1,\ldots,K\}$.
In the final step of the proof we further show that all solutions with $l_{k-1,1}\neq l_{k,1}$ for even integers $k$ may be neglected as well. Therefore we first prove that for $\delta_1,\delta_2\in\{0,1\}$ and any constant $R\in\mathbb{R}$ we have
\begin{equation}
 \label{233}
 \sum_{0<|j_{k-1,1}|,|j_{k,1}|\leq N^K}\sum_{0\leq l_{k,1}\leq l_{k-1,1}\leq 2N}|j_{k-1,1}|^{-1}|j_{k,1}|^{-1}=\mathcal{O}(N)
\end{equation}
where the sum is taken over all solutions with 
\begin{equation}
 \label{234}
 \left|(-1)^{\delta_1}j_{k-1,1}M_{l_{k-1,1},1}+(-1)^{\delta_2}j_{k,1}M_{l_{k,1},1}-R\right|<\frac{1}{2}M_{l_{k,1},1}.
\end{equation}
With $R=\alpha_{l_{k-1,1}}M_{l_{k-1,1},1}$ this condition is equivalent to
\begin{equation}
 \label{235}
 \left|\left((-1)^{\delta_1}j_{k-1,i}-\alpha_{l_{k-1,1}}\right)\frac{M_{l_{k-1,1},1}}{M_{l_{k,1},1}}+(-1)^{\delta_2}j_{k,1}\right|<\frac{1}{2}.
\end{equation}
For each choice of $l_{k-1,1}$,$l_{k,1}$ and $j_{k-1,1}$ there exists at most one choice for $j_{k,1}=j_{k,1}(j_{k-1,1})$. 
Thus we estimate the left-hand side of (\ref{233}) by
\begin{multline}
 \label{236}
 \sum_{0\leq l_{k,1}\leq l_{k-1,1}\leq 2N}\sum_{0<|j_{k-1,1}|,|j_{k,1}|\leq N^K}|j_{k-1,1}|^{-1}|j_{k,1}|^{-1}\\
 \begin{aligned}
  \leq & \sum_{0\leq l_{k,1}\leq l_{k-1,1}\leq 2N}\sum_{\substack{\lceil\alpha_{l_{k-1,1}}\rceil+1\leq (-1)^{\delta_1}j_{k-1,1}\leq N^K,\\ j_{k-1,1}\neq 0}}|j_{k-1,1}|^{-1}|j_{k,1}(j_{k-1,1})|^{-1}\\
  & + \sum_{0\leq l_{k,1}\leq l_{k-1,1}\leq 2N}\sum_{\substack{-N^K\leq (-1)^{\delta_1}j_{k-1,1}\leq\lfloor\alpha_{l_{k-1,1}}\rfloor-1,\\ j_{k-1,1}\neq 0}} |j_{k-1,1}|^{-1}|j_{k,1}(j_{k-1,1})|^{-1}\\
  & + \sum_{0\leq l_{k,1}\leq l_{k-1,1}\leq 2N}\sum_{\lfloor\alpha_{l_{k-1,1}}\rfloor\leq (-1)^{\delta_1}j_{k-1,1}\leq\lceil\alpha_{l_{k-1,1}}\rceil}|j_{k-1,1}|^{-1}|j_{k,1}(j_{k-1,1})|^{-1}.
 \end{aligned}
\end{multline}
There exists at most two $j_{k-1,1}$ such that $|(-1)^{\delta_1}j_{k-1,1}-\alpha_{l_{k-1,1}}|<1$. Observe that in this case we have $\lfloor\alpha_{l_{k-1,1},1}\rfloor\leq (-1)^{\delta_1}j_{k-1,1}\leq \lceil\alpha_{l_{k-1,1},1}\rceil$.
Therefore in any other case we have $|(-1)^{\delta_1}j_{k-1,1}-\alpha_{l_{k-1,1}}|\geq z$ for some $z\geq 1$. By (\ref{235}) we get 
\begin{equation*}
 |j_{k,1}(j_{k-1,1})|\geq \min(zq^{l_{k-1,1}-l_{k,1}}-1/2,1)\geq \min(q^{l_{k-1,1}-l_{k,1}}+z-3/2,1).
\end{equation*} 
Therefore using Cauchy-Schwarz inequality we obtain
\begin{multline}
 \label{237}
 \sum_{0\leq l_{k,1}\leq l_{k-1,1}\leq 2N}\sum_{0<|j_{k-1,1}|,|j_{k,1}|\leq N^K}|j_{k-1,1}|^{-1}|j_{k,1}|^{-1}\\
 \begin{aligned}
  \leq & \sum_{0\leq l_{k,1}\leq l_{k-1,1}\leq 2N}2\left(\sum_{(-1)^{\delta_1}j_{k,1}=\min(zq^{l_{k-1,1}-l_{k,1}}-1/2,1)}^{\infty}j_{k,1}^{-2}\right)^{1/2}\\
  & \cdot\left(\left(\sum_{(-1)^{\delta_1}j_{k-1,1}=\lceil\alpha_{l_{k-1,1}}\rceil+1}^{N^K}j_{k-1,1}^{-2}\right)^{1/2}+\left(\sum_{(-1)^{\delta_1}j_{k-1,1}=-N^K}^{\lfloor\alpha_{l_{k-1,1}}\rfloor-1}j_{k-1,1}^{-2}\right)^{1/2}\right)\\
  & +\sum_{0\leq l_{k,1}\leq l_{k-1,1}\leq 2N}\sum_{\lfloor\alpha_{l_{k-1,1}}\rfloor\leq (-1)^{\delta_1}j_{k-1,1}\leq\lceil\alpha_{l_{k-1,1}}\rceil}|j_{k-1,1}|^{-1}|j_{k,1}(j_{k-1,1})|^{-1}\\
  \leq & \sum_{0\leq l_{k,1}\leq l_{k-1,1}\leq 2N}C\max(q^{l_{k-1,1}-l_{k,1}}-3/2,1)^{-1/2}\\
  & +\sum_{0\leq l_{k,1}\leq l_{k-1,1}\leq 2N}\sum_{\lfloor\alpha_{l_{k-1,1}}\rfloor\leq (-1)^{\delta_1}j_{k-1,1}\leq\lceil\alpha_{l_{k-1,1}}\rceil}|j_{k-1,1}|^{-1}|j_{k,1}(j_{k-1,1})|^{-1}.
 \end{aligned}
\end{multline}
Since the first term is bounded by $CN$ for some constant $C>0$ which only depends on $q$ it only remains to bound the second term. Let $l_{k-1,1}$ be some fixed integer.
Now let $\tilde{l}_{k-1,1}$ be the largest integer such that (\ref{235}) holds with $\tilde{l}_{k-1,1}=l_{k,1}$ for some arbitrary $l_{k,1}$. With $|(-1)^{\delta_2}j_{k,1}|\geq 1$ we get 
\begin{equation*}
 \left|(-1)^{\delta_1}j_{k-1,1}-\alpha_{l_{k-1,1}}\right|\frac{M_{l_{k-1,1},1}}{M_{\tilde{l}_{k-1,1}}}>1/2.
\end{equation*}
For some $l_{k,1}\leq \tilde{l}_{k-1,1}-\log_q(2)-1$ we obtain
\begin{equation*}
 \left|(-1)^{\delta_1}j_{k-1,1}-\alpha_{l_{k-1,1}}\right|\frac{M_{l_{k-1,1},1}}{M_{l_{k,1},1}}>\frac{1}{2}\frac{M_{\tilde{l}_{k-1,1},1}}{M_{l_{k,1},1}}\geq q^{\tilde{l}_{k-1,1}-l_{k,1}}
\end{equation*}
and therefore $|j_{k,1}|>1/2\cdot q^{\tilde{l}_{k-1,1}-l_{k,1}}$. Using a similar argumentation as above we have
\begin{multline*}
 \sum_{0\leq l_{k,1}\leq l_{k-1,1}\leq 2N}\sum_{\lfloor\alpha_{l_{k-1,1}}\rfloor\leq (-1)^{\delta_1}j_{k-1,1}\leq\lceil\alpha_{l_{k-1,1}}\rceil}|j_{k-1,1}|^{-1}|j_{k,1}(j_{k-1,1})|^{-1}\\
 \begin{aligned}
  \leq & 2(\log_q(2)+1)(2N+1)+C\sum_{l=0}^{2N}\sum_{l'=0}^{l}q^{l'-l}\\
  \leq & CN
 \end{aligned}
\end{multline*}
for some constant $C>0$ which only depends on $q$. Plugging this into (\ref{237}) we get (\ref{233}).
With Lemma \ref{40} we get
\begin{multline}
 \label{242}
 \left|\sum_{\substack{(\delta_1,\ldots,\delta_K,j_1,\ldots,j_K,\\ n_1,\ldots,n_K)\in\mathcal{A}_{\{2,4,\ldots,K\}},\\l_{k-1,1}\neq l_{k,1} \textnormal{ for some even }k,\\ \Delta\in\mathcal{D}}}\frac{1}{2^K}\frac{1}{N^{K/2+1}}\prod_{k=1}^Ka'_{j_k}\right|\\
 \begin{aligned}
 \leq & \sum_{\substack{\delta_1,\ldots,\delta_K\in\{0,1\}}}\sum_{n_1,\ldots,n_K\in\{1,\ldots,N\}} C\\
 & \cdot \sum_{0<|j_{1,1}|,\ldots,|j_{K,1}|\leq N^K}\prod_{k=1}^K|j_{k,1}|^{-1}\mathbf{1}_{\substack{|\sum_{k=1}^H(-1)^{\delta_{\pi_1(k)}}j_{\pi_1(k),1}M_{l_{k,1},1}|<1/2\cdot M_{l_{H,1},1}\\ \textnormal{if and only if } H \textnormal{ is even}}}\\
 & \cdot \sum_{0<|j_{1,2}|,\ldots,|j_{K,2}|\leq N^K}\prod_{k=1}^K|j_{k,2}|^{-1}\mathbf{1}_{\substack{|\sum_{k=1}^H(-1)^{\delta_{\pi_2(k)}}j_{\pi_2(k),2}M_{l_{k,2},2}|<1/2\cdot M_{l_{H,2},2}\\ \textnormal{if and only if } H \textnormal{ is even}}}
 \end{aligned}
\end{multline}
for some constant $C>0$ independent of $N$. For any fixed $\delta_1,\ldots,\delta_K$, $n_1,\ldots,n_K$ and $j_{\pi_2(k),2}$ for $\pi_2(k)\leq \pi_2(H-2)$ there exists at most one $j_{\pi_2(H),2}$ for any $j_{\pi_2(H-1),2}$ and vice versa such that
\begin{equation*}
 \left|\sum_{k=1}^H(-1)^{\delta_{\pi_2(k)}}j_{\pi_2(k),2}M_{l_{k,2},2}\right|<\frac{1}{2}M_{l_{H,2},2}
\end{equation*}
holds. By Cauchy-Schwarz inequality we obtain
\begin{multline}
 \label{243}
 \sum_{0<|j_{1,2}|,\ldots,|j_{K,2}|\leq N^K}\prod_{k=1}^K|j_{k,2}|^{-1}\mathbf{1}_{\substack{|\sum_{k=1}^H(-1)^{\delta_{\pi_2(k)}}j_{\pi_2(k),2}M_{l_{k,2},2}|<1/2\cdot M_{l_{H,2},2}\\ \textnormal{if and only if } H \textnormal{ is even}}}\\
 \leq \left(2\sum_{j=1}^{\infty}j^{-2}\right)^{K/2}\leq C
\end{multline}
for some constant $C>0$ independent of $N$.
Let $H$ be the smallest even positive integer such that $l_{H-1,1}\neq l_{H,1}$. We use condition (\ref{251}) and (\ref{238}) for $k=H$ and apply (\ref{233}) with $R=\sum_{m=1}^{k-2}(-1)^{\delta_{\pi_1(m)}}j_{\pi_1(m),1}M_{l_{m,1},1}$ for any other even integer $k$. By (\ref{242}) and (\ref{243}) we conclude
\begin{equation}
 \label{244}
 \left|\sum_{\substack{(\delta_1,\ldots,\delta_K,j_1,\ldots,j_K,\\ n_1,\ldots,n_K)\in\mathcal{A}_{\{2,4,\ldots,K\}},\\l_{H-1,1}\neq l_{H,1},\\ \Delta\in\mathcal{D}}}\frac{1}{2^K}\frac{1}{N^{K/2+1}}\prod_{k=1}^Ka'_{j_k}\right|=o(1)
\end{equation}
since $K\cdot N^K\cdot q^{-2CN^{1-\varepsilon}}\leq 1/2\cdot q^{CN^{1-\varepsilon}}$ for sufficiently large $N$. Therefore we may assume $l_{k',1}=l_{k,1}$ for any $\{k,k'\}\in\Delta$.

We now show that we may also assume $l_{k',2}=l_{k,2}$ for $\{k,k'\}\in\Delta$. This immediately follows by proving $n_{k+1}=n_{k'}$ resp. $n_{k'+1}=n_{k}$ for any $\mathcal{S}=\{k,k'\}\in\Delta$ which we are going to show by induction. For $K=2$ this is trivial. For $K\geq 4$ there exists some $\{H-1,H\}\in\Delta$ since $\Delta$ is non-crossing. Furthermore by $l_{\pi_1^{-1}(H-1),1}=l_{\pi_1^{-1}(H),1}$ we have $n_{H-1}=n_{H+1}$.
Set $\tilde{k}=k$ for $k\in\{1,\ldots,H-1\}$ and $\tilde{k}=k-2$ for $k\in\{H+1,\ldots,K-1\}$. Thus it may be reduced to the case $K-2$ and the conclusion follows by induction hypothesis.
Observe that there are $N$ choices for $n_1$ and for $k\in\{2,\ldots,K\}$ the number of choices for $n_{k+1}$ is either $1$, in the case $\{k',k\}\in\Delta$ for some $k'\in\{1,\ldots,k-1\}$, or $N$ otherwise. Thus we have
\begin{equation}
 \label{249}
 \#n_1,\ldots,n_K=N^{K/2+1}.
\end{equation}
Now let $\{k,k'\}\in\Delta$. With (\ref{245}) it is easy to see that $(-1)^{\delta_{k'}}j_{k',i}+(-1)^{\delta_{k}}j_{k,i}=0$ for any $i\in\{1,2\}$ resp. $j_{k'}=\pm j_{k}$. Hence we get
\begin{multline}
 \label{246}
 a'_{j_{k'}}=\prod_{i=1}^2\frac{N^K+1-|j_{k,i}|}{N^K+1}\int_{[0,1)^2}f(x)\cos(2\pi\langle j_{k'},x\rangle)\,dx\\
 =\prod_{i=1}^2\frac{N^K+1-|j_{k,i}|}{N^K+1}\int_{[0,1)^2}f(x)\cos(2\pi\langle j_{k},x\rangle)\,dx=a'_{j_{k}}.
\end{multline}
Therefore we obtain
\begin{equation}
 \sum_{\substack{\delta_1,\ldots,\delta_K,j_1,\ldots,j_K,\\ n_1,\ldots,n_K\in\mathcal{A}_{\{2,4,\ldots,K\}},\\l_{k-1,1}=l_{k,1} \forall k\in\{2,4,\ldots,K\},\\ \Delta\in\mathcal{D}}}\frac{1}{2^K}\frac{1}{N^{K/2+1}}\prod_{k=1}^Ka'_{j_k}= |\mathcal{D}|\left(\sum_{\substack{j\in(\mathbb{Z}\backslash \{0\})^2,\\ ||j||_{\infty}\leq N^K}}a'^2_j\right)^{K/2}.
\end{equation}
By (\ref{221}) and $||f-p||_2^2\leq CN^{-K}$ we have $1\geq\sum_{j\in(\mathbb{Z}\backslash \{0\})^2}a'^2_j\geq 1-CN^{-K}$. With (\ref{225}), (\ref{247}), (\ref{239}), (\ref{240}), (\ref{241}) and (\ref{244}) we get (\ref{250}) which finally concludes the proof.

\vspace{5ex}

\section{Examples}
\label{S52}

\subsection{Superlacunary sequences}

Consider a random matrix ensemble where the first generating lacunary sequence is superlacunary, i.e. let $(M_{n,1})_{n\geq 1}$ be some sequence of positive integers such that 
\begin{equation*}
 q_n=\frac{M_{n+1,1}}{M_{n,1}}\to\infty 
\end{equation*}
for $n\to\infty$. In order to prove weak convergence to the semicircle law we have to verify (\ref{251}).
Repeating the proof of (\ref{233}) with $R=0$ we obtain by applying Cauchy-Schwarz inequality
\begin{multline*}
 \sum_{\substack{n,n'\in\{1,\ldots,2N\},\\ n>n'}}\sum_{j,j'\in\{1,\ldots,N^K\}}(jj')^{-1}\mathbf{1}_{|jM_{n,1}-j'M_{n',1}|<1/2\cdot q^{-CN^{1-\varepsilon}}M_{n',1}}\\
 \begin{aligned}
  \leq & \sum_{\substack{n,n'\in\{1,\ldots,2N\},\\ n>n'}}\sum_{j,j'\in\{1,\ldots,N^K\}}(jj')^{-1}\mathbf{1}_{|jM_{n,1}/M_{n',1}-j'|<1/2}\\
  \leq & \sum_{\substack{n,n'\in\{1,\ldots,2N\},\\ n>n'}}\left(\sum_{j=1}^{\infty}j^{-2}\right)^{1/2}\left(\sum_{j'=\max(q_{n'}^{n-n'}-1/2,1)}^{\infty}(j')^{-2}\right)^{1/2}\\
  \leq & \sum_{\substack{n,n'\in\{1,\ldots,2N\},\\ n>n'}}C\min\left((q_{n'}^{n-n'}-1/2)^{-1/2},1\right)
 \end{aligned}
\end{multline*}
for some constants $C>0$ independent of $N$. With $\tilde{n}=\min\{n:q_n\geq 2\}$ we have $(q_{n'}-3/2)^{-1/2}\leq 2(q_{n'})^{-1/2}\leq 2$ for $n'\geq \tilde{n}$. Thus we get
\begin{multline*}
 \sum_{\substack{n,n'\in\{1,\ldots,2N\},\\ n>n'}}\sum_{j,j'\in\{1,\ldots,N^K\}}(jj')^{-1}\mathbf{1}_{|jM_{n,1}-j'M_{n',1}|<1/2\cdot q^{-CN^{1-\varepsilon}}M_{n',1}}\\
 \leq C\tilde{n}+\sum_{n=1}^{2N}\sum_{k=1}^{n-1}Cq_{n-1}^{-k/2}.
\end{multline*}
Since $q_n^{-1}\to 0$ for $n\to\infty$ we easily verify (\ref{251}).

\subsection{The sequence $2^n$}

We now are going to show that besides the Hadamard gap condition (\ref{17}) further conditions on the generating lacunary sequences are necessary. Therefore we prove that for the sequences $M_{n,1}=M_{n,2}=2^n$ which do not satisfy (\ref{251}) the mean empirical eigenvalue distribution does not converge to the semicircle law
in general. Here we consider $f(x_1,x_2)=1/\sqrt{2}\cdot (\cos(2\pi(x_1+x_2))+\cos(4\pi(x_1+x_2)))$. We restrict ourselves to the case $K=4$ and repeat the proof of Theorem \ref{253} until (\ref{241}). Observe that there are precisely two non-crossing pair partitions of $\{1,2,3,4\}$ which are $\{\{1,2\},\{3,4\}\}$ resp. $\{\{1,4\},\{2,3\}\}$. Thus $\mathbb{E}[1/N\cdot X_N^4]$ is the sum over solutions $\delta_1,\ldots,\delta_4$,$j_1,\ldots,j_4$ and $n_1,\ldots,n_4$ such that
\begin{equation*}
 (-1)^{\delta_1}j_{1,1}M_{n_1+n_2,1}+(-1)^{\delta_2}j_{2,1}M_{n_2+n_3,1}=(-1)^{\delta_3}j_{3,1}M_{n_3+n_4,1}+(-1)^{\delta_4}j_{4,1}M_{n_4+n_1,1}=0
\end{equation*}
or
\begin{equation*}
 (-1)^{\delta_1}j_{1,1}M_{n_1+n_2,1}+(-1)^{\delta_4}j_{4,1}M_{n_4+n_1,1}=(-1)^{\delta_2}j_{2,1}M_{n_2+n_3,1}+(-1)^{\delta_3}j_{3,1}M_{n_3+n_4,1}=0.
\end{equation*}
Without loss of generality by (\ref{238}) we may assume that both cases are distinct. Hereafter we only focus on the first case. The second case is similar.
We further get
\begin{eqnarray*}
 (j_{2,1},j_{2,2}) & = & (-1)^{\delta_1+\delta_2}(2^{n_1-n_3}j_{1,1},2^{|n_1-n_2|-|n_2-n_3|}j_{1,2}),\\
 (j_{4,1},j_{4,2}) & = & (-1)^{\delta_3+\delta_4}(2^{n_3-n_1}j_{3,1},2^{|n_3-n_4|-|n_4-n_1|}j_{3,2}).
\end{eqnarray*}
Therefore we have
\begin{equation*}
 \begin{aligned}
  \lim_{N\to\infty}\frac{1}{N}\mathbb{E}[X_N^4] \:\:\: & = \:\:\: \frac{2}{N^3}\sum_{n_1,n_3=1}^N\left(\sum_{n=1}^N\mathbb{E}\left[f\left(2^{n_3-n_1\vee 0}x_1,2^{|n-n_3|-|n_1-n|\vee 0}x_2\right)\right.\right.\\
  & \quad\quad\quad\quad\quad\quad\quad\quad\quad\quad\left.\left.\cdot f\left(2^{n_1-n_3\vee 0}x_1,2^{|n_1-n|-|n-n_3|\vee 0}x_2\right)\right]\right)^2.
 \end{aligned}
\end{equation*}
For $n_1=n_3$ we easily observe
\begin{equation*}
 \frac{2}{N^3}\sum_{n_1=1}^N\left(\sum_{n=1}^N\mathbb{E}[f^2]\right)^2=2=C_2.
\end{equation*}
Thus it is enough to show
\begin{equation}
 \label{254}
 \sum_{n=1}^N\mathbb{E}\left[f\left(2^{n_3-n_1\vee 0}x_1,2^{|n-n_3|-|n_1-n|\vee 0}x_2\right)f\left(2^{n_1-n_3\vee 0}x_1,2^{|n_1-n|-|n-n_3|\vee 0}x_2\right)\right]\neq 0
\end{equation}
for some $n_1\neq n_3$ in order to prove that the mean empirical eigenvalue distribution does not converge to the semicircle law.
Now choose some $n_1=n_3+1>1$. Then we have
\begin{equation*}
 |n_1-n|-|n-n_3|=
 \begin{cases}
  1, & n\leq n_3,\\
  -1, & n\geq n_1.
 \end{cases}
\end{equation*}
Plugging this into the left-hand side of (\ref{254}) and using the definition of $f$ yields
\begin{equation*}
 \sum_{n=1}^{n_3}\mathbb{E}\left[f(x_1,x_2)f\left(2x_1,2x_2\right)\right]+\sum_{n=n_1}^N\mathbb{E}\left[f(x_1,2x_2)f\left(2x_1,x_2\right)\right]=\frac{n_3}{2}>0.
\end{equation*}
Therefore we obtain
\begin{equation*}
 \lim_{N\to\infty}\frac{1}{N}\mathbb{E}[X_N^4]\neq 2.
\end{equation*}
By assumption on $f$ we trivially have
\begin{equation*}
 \lim_{N\to\infty}\frac{1}{N}\mathbb{E}[X_N^2]=\mathbb{E}[f^2]=1.
\end{equation*}
We conclude that the limiting mean empirical eigenvalue distribution does not have the density $1/R\pi\cdot\sqrt{R^2-x^2}\cdot\mathbf{1}_{x^2\leq R^2}$ for any $R>0$ and hence it is not a semicircle law.

\vspace{3ex}

\small{DEPT. OF MATHEMATICS, BIELEFELD UNIV., P.O.Box 100131, 33501 Bielefeld, Germany\\
\textit{E-Mail address:} \url{tloebbe@math.uni-bielefeld.de}}

\end{document}